\documentclass[11pt]{amsart}
\usepackage{amsmath,amsthm,amsfonts}
\usepackage{graphicx}

\theoremstyle{definition}

\newcommand{\gp}[1]{{\left\langle #1 \right\rangle}}

\title[A new key exchange protocol]{A new key exchange protocol based on the decomposition problem}

\date{}

\author{Vladimir Shpilrain}
\address{Department of Mathematics, The City  College  of New York, New York,
NY 10031} \email{shpil@groups.sci.ccny.cuny.edu}

\author{Alexander Ushakov}
\address{Department of Mathematics, Stevens Institute of Technology, Hoboken, NJ 07030}
\email{aushakov@mail.ru}
\thanks{Research of the first author was partially supported by
the NSF grant DMS-0405105.}

\begin{document}

\maketitle

\begin{abstract}
In this paper we present   a new key establishment protocol based on
the decomposition problem in non-commutative groups which is: given
two elements $w, w_1$ of the platform group  $G$ and two
 subgroups $A, B \subseteq G$ (not necessarily distinct), find  elements
 $a \in A, ~b \in B$ such that $w_1 = a w b$. Here we introduce two new ideas that improve the security
 of key establishment protocols based on the decomposition problem.
 In particular, we conceal (i.e., do not publish explicitly) one of the subgroups $A,
 B$, thus introducing an additional computationally hard problem
 for the adversary, namely, finding the centralizer of a given
 finitely generated subgroup.

\end{abstract}

\section{Introduction}

  In search of a more efficient and/or secure alternative to established
cryptographic protocols (such as RSA), several authors have come up
with public key establishment protocols as well as with complete
public key cryptosystems based on allegedly hard {\it search
problems} from combinatorial (semi)group theory, including the
conjugacy  search problem \cite{AAG, KLCHKP}, the homomorphism
search problem \cite{Grigoriev}, \cite{SZ1},  the decomposition
 search problem \cite{CKLHC, KLCHKP, SU}, the subgroup membership search problem
\cite{SZ2}.

 In this paper, we focus on the decomposition search problem which we
 subsequently call just the decomposition problem. The problem is:
 given two elements $w, w_1$ of the platform group  $G$ and two
 subgroups $A, B \subseteq G$ (not necessarily distinct), find  elements
 $a \in A, ~b \in B$ such that $w_1 = a w b$.

 It is straightforward to arrange a key establishment protocol based on
 this problem (see \cite{CKLHC, KLCHKP, SU}), assuming that $ab=ba$ for any
 $a \in A, ~b \in B$:

\medskip\noindent{\bf (0)} One of the parties (say, Alice) publishes
a random element $w \in G$ (the ``base" element).

\medskip\noindent{\bf (1)} Alice chooses   $a_1, a_2 \in A$
(Alice's private keys) and sends $a_1 w a_2$ to Bob.

\medskip\noindent{\bf (2)} Bob chooses $b_1, b_2 \in B$
(Bob's private keys) and   sends $b_1 w b_2$ to Alice.

\medskip\noindent{\bf (3)} Alice computes
$$K_a = a_1 b_1 w b a_2 b_2$$
and Bob computes
$$K_b = b_1 a_1 w a_2 b_2.$$
If $a_ib_i=b_ia_i$,  then  $K_a=K_b$ in $G$. Thus Alice and Bob have
a {\em shared secret key}.
\medskip

Security of such a protocol
  will, of course, depend on a particular platform group  $G$ (at
  the very least, $G$ has to be non-commutative). It appears
that   for braid groups (which are a popular choice for the
platform), the so-called {\it length attacks} present a serious
threat, see e.g. \cite{GKTTV2, HS, HT, MSU}.

 In this paper, we introduce two new ideas that improve the security
 of key establishment protocols based on the decomposition problem:

\medskip\noindent{\bf (i)}  We conceal one of the subgroups $A, B$.

\medskip\noindent{\bf (ii)} We make Alice   choose her left private key $a_1$ from
one of the subgroups $A, B$, and her right private key $a_2$  from
the other subgroup. Same for Bob.
\medskip

These two improvements together
 will obviously foil any length attacks. We give a complete
 description of our protocol in the following Section
 \ref{se:protocol}; here we just sketch the main idea.

Let $G$ be a group and $g \in G$. Denote by $C_G(g)$ the {\em
centralizer} of $g$ in $G$, i.e., the set of elements
 $h \in G$ such that $hg = gh$.
 For $S = \{g_1,\ldots,g_k \} \subseteq G$,   $C_G(g_1,\ldots,g_k)$   denotes
the centralizer of $S$ in $G$, which is the intersection of  the
centralizers $C_G(g_i), i=1,...,k$.

 Now, given a public $w \in G$, Alice
privately selects $a_1 \in G$ and publishes a subgroup   $B
\subseteq C_G(a_1)$ (we explain why computing $B$ is easy).
Similarly, Bob privately selects $b_2 \in G$ and publishes a
subgroup $A \subseteq C_G(b_2)$. Alice then selects $a_2 \in A$ and
sends $w_1=a_1wa_2$ to Bob, while Bob selects $b_1 \in B$ and sends
$w_2=b_1wb_2$ to Alice.

Thus,  in the first  transmission, say, the adversary faces the
problem of finding $a_1, a_2$ such that $w_1 = a_1 w a_2$, where
$a_2 \in A$, but there is no explicit indication of where to choose
$a_1$ from. Therefore, before arranging something like a length
attack in this case, the adversary would have to compute the
centralizer $C_G(B)$ first (because $a_1 \in C_G(B)$), which is
usually a hard problem by itself.

\section{The protocol}
\label{se:protocol}

In this section we give a formal description of our protocol, but
first we introduce one more piece of notation. As it is common in
public key exchange based on abstract groups, when transmitting an
element $g \in G$ of a group, one actually uses its {\it normal
form} $N(g)$ which is a sequence of symbols uniquely defined for a
given $g$. A specific way of constructing such a sequence depends,
of course, on a particular platform group $G$ which we discuss in
subsequent sections of our paper.

Our protocol is the  following sequence of steps.

\medskip \emph{\textbf{Protocol:}}

\begin{enumerate}

    \item[\bf(1)]
Alice chooses an element $a_1 \in G$ of length $l$, chooses a
subgroup of $C_G(a_1)$, and publishes its generators $A = \{
\alpha_1,\ldots,\alpha_k\}$ (see the following subsection
\ref{parameters} for specifications).
    \item[\bf(2)]
Bob chooses an element $b_2 \in G$ of length $l$, chooses a subgroup
of $C_G(b_2)$, and publishes  its generators $B =
\{\beta_1,\ldots,\beta_m \}$ (see the following subsection
\ref{parameters} for specifications).
    \item[\bf(3)]
Alice chooses a random element $a_2$ from
$\gp{\beta_1,\ldots,\beta_m}$ and sends the normal form $P_A = N(a_1
w a_2)$ to Bob.
    \item[\bf(4)]
Bob chooses a random element $b_1$ from
$\gp{\alpha_1,\ldots,\alpha_k}$ and sends  the normal form  $P_B =
N(b_1 w b_2)$ to Alice.
    \item[\bf(5)]
Alice computes $K_A = a_1 P_B a_2$.
    \item[\bf(6)]
Bob computes $K_B = b_1 P_A b_2$.
\end{enumerate}

Since $a_1b_1 = b_1a_1$ and $a_2b_2 = b_2a_2$, we have $K = K_A =
K_B$, the shared secret  key.

\subsection{Suggested values of parameters}
\label{parameters}

We suggest to use the following values of parameters in the above
protocol: $G=B_n$, the group of braids on $n$ strands (see our
Section \ref{braids});  $n = 64$; $l = 1024$. At Step (1) of the
protocol    Alice   generates $(a_1, A)$   and at Step (2) Bob
generates $(b_2, B)$, both using the algorithm from \cite{FGM} for
computing centralizers (actually, there is no need to  compute the
whole centralizer, just a couple of elements are   sufficient).

\section{Requirements on the platform  group $G$}
\label{requirements}

In this section we discuss possible attacks on
 the protocol described in the previous section,
 and also put together  some requirements on the platform  group $G$.

To break the protocol it is sufficient to find either Alice's or
Bob's private key which may be accomplished as follows:
\begin{enumerate}
    \item[]
{\bf Attack on Alice's private key.} Find an element $a_1'$ which
commutes with every element of the subgroup $\gp{A}$ and an element
$a_2' \in \gp{B}$, such that $P_A = N(a_1' w a_2')$. The pair
$(a_1',a_2')$ is equivalent to $(a_1,a_2)$. (That means, $a_1' w
a_2'=a_1 w a_2$, and therefore the pair $(a_1',a_2')$ can be used by
the adversary to get the shared secret  key.)
    \item[]
{\bf Attack on Bob's private key.} Find an element $b_1' \in \gp{A}$
and an element $b_2'$ which commutes with every element of  the
subgroup $\gp{B}$, such that $P_B = N(b_1' w b_2')$. The pair
$(b_1',b_2')$ is equivalent to $(b_1,b_2)$.
\end{enumerate}

Consider the attack on  Alice's private key (the other one is
similar). The most obvious way to carry out such an attack is the
following:
\begin{enumerate}
    \item[{\bf (A1)}]
Compute the centralizer $C_G(A)$.
    \item[{\bf (A2)}]
Solve the search version of the membership problem in the double
coset $C_G(A) \cdot w \cdot \gp{B}$
\end{enumerate}

To make the protocol secure, we want both  these problems to be
computationally hard. For the problem (A2) to be hard, it is
necessary  for the centralizer $C_G(A)$ to  be large. Otherwise, the
adversary can use the ``brute force" attack, i.e.,  enumerate all
elements of $C_G(A)$ and find candidates for $b_2'$ (assuming that
the decisional membership problem in the subgroup $B$ is efficiently
solvable).

Thus the platform group $G$ should satisfy   at least the following
properties in order for our key establishment protocol  to be
efficient and secure.

\begin{enumerate}
    \item[\bf(P1)]
$G$ should be a non-commutative group of exponential growth. The
latter means that the number of elements of length $n$ in $G$ is
exponential in $n$; this is needed to prevent attacks by complete
exhaustion of the key space.
    \item[\bf(P2)]
There should  be an efficiently computable normal form for elements
of $G$.
    \item[\bf(P3)]
It should  be computationally easy to perform group operations
(multiplication and inversion) on normal forms.
    \item[\bf(P4)]
It should  be computationally easy to generate pairs $(a, ~\{a_1,
\ldots, a_k\})$ such that $a a_i = a_i a$ for each $i = 1,\ldots,k$.
(Clearly, in this case the subgroup generated by $a_1, \ldots, a_k$
centralizes $a$).
    \item[\bf(P5)]
For a generic set $\{g_1,\ldots,g_k\}$ of  elements of $G$  it
should be difficult to compute
$$C(g_1,\ldots,g_n) = C(g_1) \cap \ldots \cap C(g_k).$$
    \item[\bf(P6)]
Even if $H = C(g_1,\ldots,g_n)$ is computed, it should  be hard to
find $x \in H$ and $y \in H_1$ (where $H_1$ is some fixed subgroup
given by a generating set) such that $x w y = w'$, i.e., to solve
the membership search problem for a double coset.
\end{enumerate}

\section{Braid groups}
\label{braids}

 In this section we consider a particular class of groups,
namely braid groups, which were a popular choice for the platform of
various cryptographic protocols in the last 6-7 years, starting with
the seminal paper \cite{AAG}.

Let $B_{n}$ be the group of braids on $n$ strands and $X_n =
\{x_1,\ldots,x_{n-1}\}$  the set of standard generators. Thus,

\vskip -0.5cm

$$B_{n} = \langle x_1,\dots ,x_{n-1} ; ~x_ix_{i+1}x_i=x_{i+1}x_ix_{i+1}, ~x_ix_j=x_jx_i ~\mbox{for} ~|i-j|>1\rangle.$$

For more information on braid groups, we refer   to the monographs
\cite{Birman}, \cite{Epstein}; here we address the properties
(P1)-(P6) from the previous section.

\begin{enumerate}
    \item[(P1)]
Braid groups $B_n$ are non-commutative groups of exponential growth
if $n\ge 3$.

    \item[(P2)]
There are several known normal forms for  elements of $B_n$,
including Garside normal form (see \cite{Birman}) and Birman-Ko-Lee
normal form \cite{BKL}. Both of these forms are efficiently
computable (in quadratic time with respect to the length of a given
element).

    \item[(P3)]
There are quadratic time algorithms to multiply or invert normal
forms of elements of $B_n$.

    \item[(P4)]
It is not so easy to compute the whole centralizer of an element $g$
of $G$ (cf. \cite{GMW}). The number of steps required to compute
$C_G(g)$ is proportional to $|SSS(g)|$, the size of the ``super
summit set" of $g$,  which is typically huge. Nevertheless, there
are   approaches to finding ``large parts" of $C_G(g)$, e.g. one can
generate a sufficiently large  part of $SSS(g)$ and pick several
elements from there, see \cite{GMW} for more details.

    \item[(P5)]
For a generic subgroup $A$ it is hard to compute $C_G(A)$. The
complexity of such computation is proportional to $|SS(A)|$, the
size of the   summit set of $A$ (see \cite{FGM}), which is typically
huge.

    \item[(P6)]
There is no known solution to the membership search problem for
double cosets $H \cdot w \cdot H'$ in braid groups. This problem, in
theory, appears to be much more complicated (for generic subgroups
$H$ and $H'$) than the conjugacy search problem.
\end{enumerate}

\section{Semantic security}
\label{security}

In this section, we discuss {\it semantic security} of a
cryptosystem that would be based on a shared key obtained in our
protocol. Semantic security   is the standard notion of security for
encryption protocols, see \cite{GoMi}.

Security of the protocol described in our Section \ref{se:protocol}
is based  on the assumption that the following problem is
computationally hard:

\begin{quote}
Given the public information $w$, $P_A$, and $P_B$ it is hard to
compute the shared key $K$.
\end{quote}

This assumption is the  {\em computational assumption of the
protocol}. The stronger {\em decisional} version of this assumption
would be:

\begin{quote}
Given  $w$, $P_A$,  and  $P_B$,   it is hard to distinguish the
shared key $K$ from a random  element of the form $awb$.
\end{quote}

 We should point out that without this decisional assumption, it may
 still be possible to design a semantically secure encryption protocol in the ``random
 oracle model" the same way it was done in \cite[Section
 3.3]{KLCHKP}, namely, by employing a hash function $H : B_n \to \{0,
 1\}^k$ from the braid group to the message space. Still, it would
 be quite interesting to find out whether or not the shared key $K$
 obtained in our key establishment protocol can be directly used for semantically
 secure encryption.

The decisional assumption above appears to be wrong   for most
choices of $w, P_A$ and  $P_B$ because of the following
consideration. Since $P_A = a_1 w a_2$, we have $a_1 = P_A a_2^{-1}
w^{-1}$. Therefore, $K = a_1 b_1 w b_2 a_2 = P_A a_2^{-1} (w^{-1}
P_B) a_2$. Hence, $K$ is a product of a public element $P_A$ and a
public element $w^{-1} P_B$ conjugated by an element from a subgroup
$\{\beta_1,\ldots,\beta_k\}$.

It seems plausible that, for some choices of the keys, elements of
this type can be distinguished from random elements of the form
$awb$ along the same lines it was done in \cite{GM} (in a different,
but similar context). Indeed, if $w^{-1} P_B$ is not a {\it pure
braid},
 then it projects to a non-trivial permutation, call it $\rho_B$,
under the natural homomorphism $\pi$ from the braid group $B_{n}$
onto the symmetric group $S_{n}$. Then the conjugate permutation
$\pi(a_2)^{-1}\rho_B\pi(a_2)$ has the same cyclic structure as
$\rho_B$   does, and this gives away some information about the
permutation $\pi(K)= \pi(P_A) \pi(a_2)^{-1}\rho_B\pi(a_2)$; for
example, from knowing $\pi(P_A)$ and the cyclic structure of
$\pi(a_2)^{-1}\rho_B\pi(a_2)$, one can get information about
possible order of the permutation $\pi(K)$.

 If both $P_A$ and $w^{-1} P_B$  are pure braids, then it is possible to use
   other homomorphisms  (e.g. pulling out a strand) to obtain some
   partial information; see \cite{GM} for details. If   $w^{-1} P_B$
   is a pure braid but $P_A$ is not, then, again, the homomorphism $\pi$
 reveals partial information about the shared key $K$.

\end{document}